# TRAJECTORY GENERATION AND DISPLAY FOR FREE FLIGHT

**Mohammad Shahzad**[*], **Félix Mora-Camino**[*,✈], **Jules G. Slama**[#] **and Karim Achaibou**[*]

\* LAAS du CNRS, Toulouse, France.
\# Programa de Engenharia Mecânica, COPPE/UFRJ, Rio de Janeiro, Brazil.
✈ Air Transportation Department, ENAC, Toulouse, France.
e-mail: shahzad@laas.fr, mora@laas.fr, jules@rionet.com.br & achaibou@laas.fr

## Abstract

*In this study a new approach is proposed for the generation of aircraft trajectories. The relative guidance of an aircraft, which is aimed to join in minimum time the track of a leader aircraft, is particularly considered. In a first place, a minimum time relative convergence problem is considered and optimal trajectories are characterized. Then the synthesis of a neural approximator for optimal trajectories is discussed. Trained neural networks are used in an adaptive manner to generate intent trajectories during operation. Finally simulation results involving two wide body aircraft are presented.*

## 1 Introduction

In order to absorb the increasing air traffic flows, the innovative concept of Free Flight ("Safe and efficient flight operations in which the aircraft have the freedom to select their path and speed in real time") has been under study over the last years [1]. The implementation of Free Flight has been made possible by the emergence of new navigation technologies such as GPS, ADS-B, TCAS and of new onboard computation capabilities [2]. Hence the absolute position of an aircraft as well as its short term intents can be communicated to other aircraft in the neighborhood through data links, while precise relative positions can be computed on line. Then aircraft should be able to realize relative maneuvers such as minimum separation crossings, mergings and meterings along common tracks. In this communication, the case of the convergence maneuver is considered more particularly.

Time and cost optimization of aircraft trajectories has been of great interest for many decades and various numerical solution techniques have been developed [3-6]. However these techniques are not in general compatible with an on-line operation which is here a necessity since the leading aircraft may modify at any time its guidance parameters (speed, heading and flight level) in accordance with new atmospheric conditions (wind and temperature) or following instructions issued by ATC. On the other hand, simple proportional navigation techniques, developed for missile homing, provide direct solutions to a similar problem [7-8]. However air transportation regulations (load factor limitations and standard maneuvers) as well as economical, structural and comfort considerations, prevent their utilization in the case considered in this study [9]. The control strategy proposed in this communication is adaptive and makes use of a neural network structure to get an on line approximation of the optimal trajectory associated to the current relative situation [10].

In section 2, the minimum time optimization problem considered is displayed and analyzed. It is shown in section 3 that the adoption of a simplified flight mechanics model for its formulation allows the analysis and subsequent characterization of the optimal trajectory in terms of turns and straight line segments.

Then in section 4 the synthesis of the neural approximator of optimal trajectories is discussed. A selected set of initial situations to cover, through neural generalization, a region corresponding to 10 minutes of flight.

In section 5, simulation results of the proposed approach for different situations are displayed.

## 2 Problem Formulation and Conditions of Optimality

In the case considered, a pursuer aircraft (*P*) is attempting to follow a leader aircraft (*L*). Both aircraft are supposed to keep their speed constant until the completion of the maneuver which is characterized by the final conditions :
- the velocity vector of the pursuer aircraft is parallel to the leader aircraft's velocity,
- the pursuer has reached the same route as the leader,
- the pursuer is at a given distance $D \geq d_{min}$ from the leader.

The geometry of the convergence is shown in Fig. 1 and the equations of motion, in relative polar co-ordinates, are [9] :

$$\dot{d} = V_L \cos(\theta - \psi_L) - V_P \cos(\theta - \psi_P) \quad (1)$$

$$\dot{\theta} = \{-V_L \sin(\theta - \psi_L) + V_P \sin(\theta - \psi_P)\}/d \quad (2)$$

$$\dot{\psi}_P = r_P \quad (3)$$

where $V_L$ and $V_P$ are the speeds of the leader and of the pursuer aircraft respectively, $d(t)$ is the instantaneous separation between the two aircraft, $\psi_L$ and $\psi_P$ are their headings, $\theta$ is the line of sight angle. $\theta$, $\psi_L$, and $\psi_P$ are measured with respect to a common earth reference.

The minimum time convergence trajectory is then solution of the optimization problem :

minimize $\int_0^{t_f} dt$ with (1), (2), (3) with the constraints :

$$\phi_{min} \leq \phi_P \leq \phi_{max} \quad (4)$$

$$d_{min} - d \leq 0 \quad (5)$$

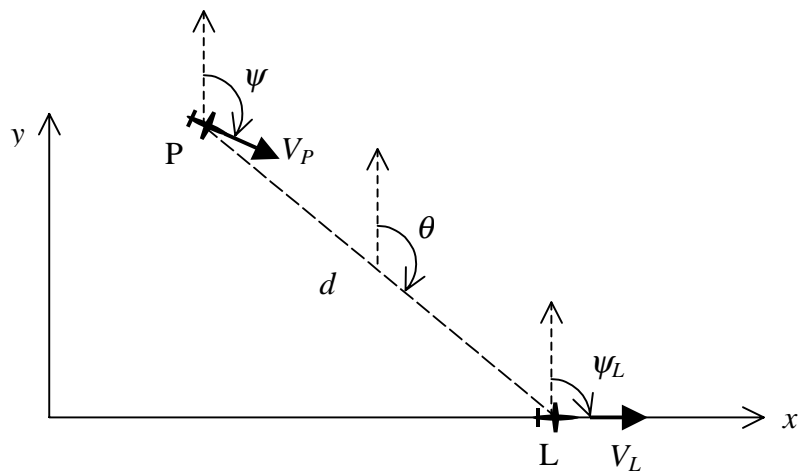

**Figure 1. Representation of relative positions.**

the initial conditions :

$$d(0) = d_0 \text{ with } d_0 \geq d_{min}$$
$$\theta(0) = \theta_0, \ \psi_P(0) = \psi_{P0} \quad (6)$$

and the final conditions of convergence :

$$\psi_P(t_f) = \psi_L, \ \theta(t_f) = \psi_L$$
$$d(t_f) = D \quad (7)$$

Here the time t is the independent variable and the control variable is the instantaneous rate of change of heading the Pursuer, $r_P$, which is associated to the bank angle, $\phi$, during a zero side-slip turn by the following relation :

$$r_P = \frac{g}{V_P} \cdot \tan\phi \quad (8)$$

The Hamiltonian for this dynamic optimisation problem is given by [11] :

$$H(\theta, d, \psi_P, \lambda_\theta, \lambda_d, \lambda_\psi, r_P) =$$
$$1 + \lambda_d \{V_L \cos(\theta - \psi_L) - V_P \cos(\theta - \psi_P)\}$$
$$+ \lambda_\theta \{-V_L \sin(\theta - \psi_L) + V_P \sin(\theta - \psi_P)\}/d$$
$$+ \lambda_\psi r_P + \mu(d_{min} - d)$$
$$(9)$$

where $\lambda_d$, $\lambda_\theta$ et $\lambda_\psi$ are adjoint variables associated to the state variables d, $\theta$ et $\psi$. In order to take into account the constraint of minimum separation an influence variable, $\mu$, is introduced:

$$\mu \geq 0 \quad \text{if} \quad (d_{min} - d) = 0 \quad (10)$$
$$\mu = 0 \quad \text{if} \quad (d_{min} - d) < 0$$

The necessary conditions of optimality are defined by Euler-Lagrange equations :

$$\dot{\lambda}_d = \frac{\lambda_\theta}{d^2} \{-V_L \sin(\theta - \psi_L) + V_P \sin(\theta - \psi_P)\}$$
$$+ \mu$$
$$(11)$$

$$\dot{\lambda}_\theta = -\frac{\lambda_\theta}{d} \{-V_L \cos(\theta - \psi_L) + V_P \cos(\theta - \psi_P)\}$$
$$- \lambda_d \{-V_L \sin(\theta - \psi_L) + V_P \sin(\theta - \psi_P)\}$$
$$(12)$$

$$\dot{\lambda}_\psi = \frac{\lambda_\theta}{d} \{V_P \cos(\theta - \psi_P)\} + \lambda_d V_P \sin(\theta - \psi_P)$$
$$(13)$$

with the final conditions of adjoint variables :

$$\lambda_\theta(t_f) = \nu_\theta, \ \lambda_\psi(t_f) = \nu_\psi, \ \lambda_d(t_f) = \nu_d \quad (14)$$

where $\nu_\theta$, $\nu_\psi$ and $\nu_d$ are the constants which must be specified in order to satisfy final conditions of the state variables.

The above system of equations should be solved numerically in order to get desired minimum time trajectories. Anyhow, by interpreting the optimality conditions, $H$ can be written as :

$$H = H(\theta, d, \psi_P, \lambda_\theta, \lambda_d) + \lambda_\psi r_P \quad (15)$$

the minimization of $H$ leads to the conditions :

$$\text{if } \lambda_\psi^* > 0 \text{ then } \phi_P^* = \phi_{min},$$
$$\text{if } \lambda_P^* < 0 \text{ then } \phi_P^* = \phi_{max}. \quad (16)$$

The case where $\lambda_\psi^* = 0$ is more complex to analyse, however this singular situation can be clarified through the manipulation of the optimality conditions [11], leading to the conclusion that when $\lambda_\psi^* = 0$ then $\phi_P^* = 0$ and $\psi_P$ remains constant.

## 3  Characterisation of Minimum Time Convergence Trajectories

The optimality conditions obtained in the previous section show that a constant speed time-optimal convergence trajectory consists of a sequence of maximum bank angle turns (right or left) linked by straight line segments.

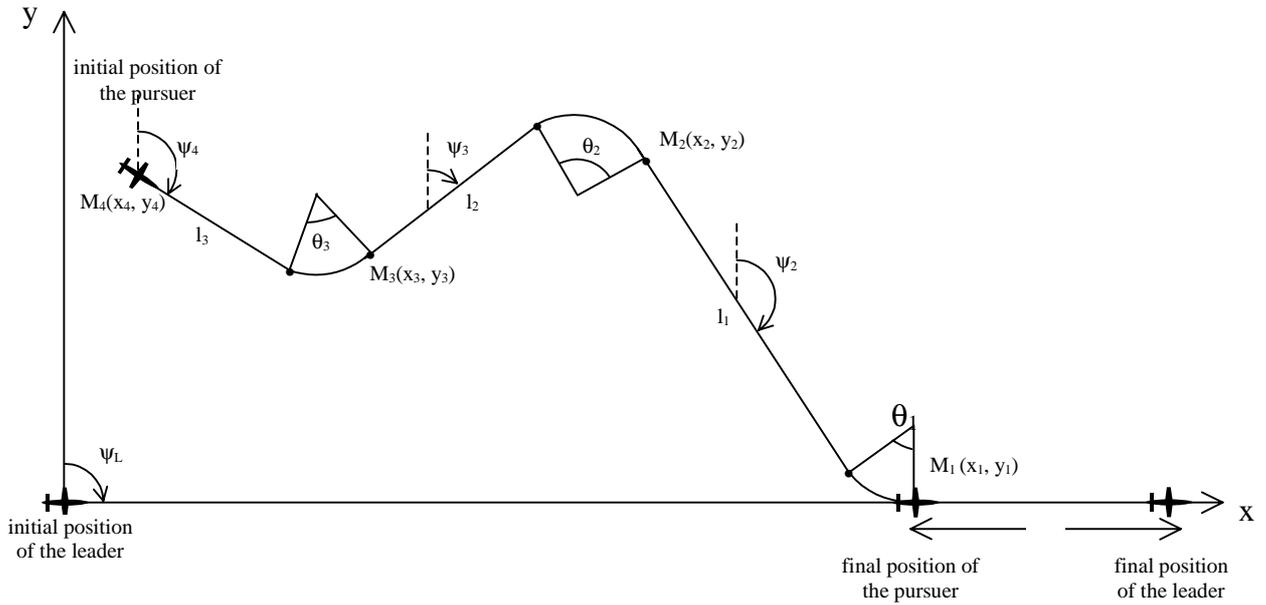

**Figure 2. Parametric representation of a regular trajectory.**

A trajectory composed of a sequence of pairs of straight line segments and maximum bank angle turns is here called "regular trajectory" and is said to be of order n if it contains (n – 1) such pairs. An example of a regular trajectory of order 4 is represented in figure 2. This class of trajectory can be parameterized by a sequence of triplets ($\varepsilon_i$, $\theta_i$, $l_i$), i = 1 to n – 1, where :

- $\varepsilon_i$ gives the direction of the turn ($\varepsilon_i$ = +1 : left turn, $\varepsilon_i$ = -1 : right turn, according to trigonometric orientation),
- $\theta_i$ is the absolute value of the turn angle realised at a nominal value of the bank angle,
- $l_i$ is the length of i$^{th}$ straight line segment.

A minimum time regular trajectory satisfying initial constraints, convergence constraints as well as minimum separation constraints is solution of the optimization problem :

$$\min_{(\varepsilon_i, \theta_i, l_i)_{i=1 \text{ à } n-1}} \sum_{k=1}^{n-1} (l_k + R_{\min}\theta_k) \text{ with } l_i \geq 0,$$

$$\theta_i \geq 0, \; \varepsilon_i = \pm 1 \quad i = 1 \text{ à n-1} \qquad (17)$$

with the following constraints :
- initial state : $x_n = x_0^P$, $y_n = y_0^P$

  and $\psi_n = \psi_0^P$ \qquad (18)

- final convergence :

$$\begin{cases} x_1 = x_L(0) + \alpha \sum_{k=1}^{n-1}(l_k + R_{\min}\theta_k) - D \\ y_1 = 0 \\ \psi_1 = \psi_L \end{cases}$$

\qquad (19)

with $\alpha = V_L/V_P$ and $x_L(0) = 0$

- minimum separation : $d \geq d_{\min}$ \qquad (20)

This last constraint is expressed according to the value of the n$^{th}$ minimum separation time $T_n^*$ given by :

$$T_n^* = -\frac{(x_n + (l_n + R_{\min}\theta_n)\sin\psi_n - x_L^{n-1})(V_P\sin\psi_n - V_L) + (y_n + (l_n + R_{\min}\theta_n)\cos\psi_n)}{(V_P\sin\psi_n - V_L)^2 + V_P^2\cos^2\psi_n} \quad (21)$$

1) if $T_n^* \leq 0$ then $\sqrt{(x_n + (l_n + R_{\min}\theta_n)\sin\psi_n - x_L^{n-1})^2 + (y_n + (l_n + R_{\min}\theta_n)\cos\psi_n)^2} \geq d_{\min}$ (22)

2) if $T_n^* \geq l_n/V_S$ then $\sqrt{(x_n - x_L^n)^2 + (y_n)^2} \geq d_{\min}$ (23)

3) if $0 < T_n^* < l_n/V_S$ then
$$(x_n + (l_n + R_{\min}\theta_n)\sin\psi_n - x_L^{n-1})^2 + (y_n + (l_n + R_{\min}\theta_n)\cos\psi_n)^2$$
$$-\frac{(x_n + (l_n + R_{\min}\theta_n)\sin\psi_n - x_L^{n-1})(V_P\sin\psi_n - V_L) + (y_n + (l_n + R_{\min}\theta_n)\cos\psi_n)}{(V_P\sin\psi_n - V_L)^2 + V_P^2\cos^2\psi_n} \geq d_{\min}^2 \quad (24)$$

From the point of view of Mathematical Programming this problem presents major difficulties :
- it is a mixed variable programming problem ($\varepsilon_i$ are binary variables while angles $\theta_i$ and lengths $l_i$ are real positive variables),
- the admissible domain generated by the different constraints is non convex and its complexity grows explosively with an increasing value of n,
- the minimum separation constraints depend on logical conditions.

The only exact method which seems applicable here is Dynamic Programming [12] and only in the forward direction (initial conditions to final convergence) because of constraints of minimum separation which destroy the property of separability and sequenciability of problem. Methods of stochastic optimisation like genetic algorithms [13] can also be considered to solve the problem numerically. In the case where an on line and on board resolution of this problem is pursued, it is clear that such numerical solution approaches cannot be considered.

It appears then interesting to investigate the techniques of Artificial Intelligence such as :
- Expert Systems [14] which can make use of a base of rules dedicated to the generation of time efficient and admissible convergence trajectories through guidance directives.
- Neural Networks which can memorise the parameters of minimum time convergence trajectories and produce from them, interpolated trajectories.

This last approach is adopted here.

## 4 Reverse Dynamic Programming and Neural Restitution of Convergence Trajectories

Since the Leader may modify at any time its guidance parameters (speed, heading and flight level) in accordance with new atmospheric conditions (wind and temperature) or following instructions issued by the traffic control service, on line trajectory generation is required here. Considering the complexity of the problem, it appears impracticable to get optimal trajectories in real time. However, the nature of the convergence maneuver imply the need to dispose at any time of guidance directives which can be immediately displayed to the pilot or taken into account by the autopilot. A neural network fulfils this requirement as it can be trained to memorize input-output relationships for a system and afterwards can be used as an interpolator to generate outputs corresponding to current inputs. In fact since the early 1990s, neural networks have been used in a variety of applications and there has been a growing interest in using neural networks to solve optimal guidance problems [15].

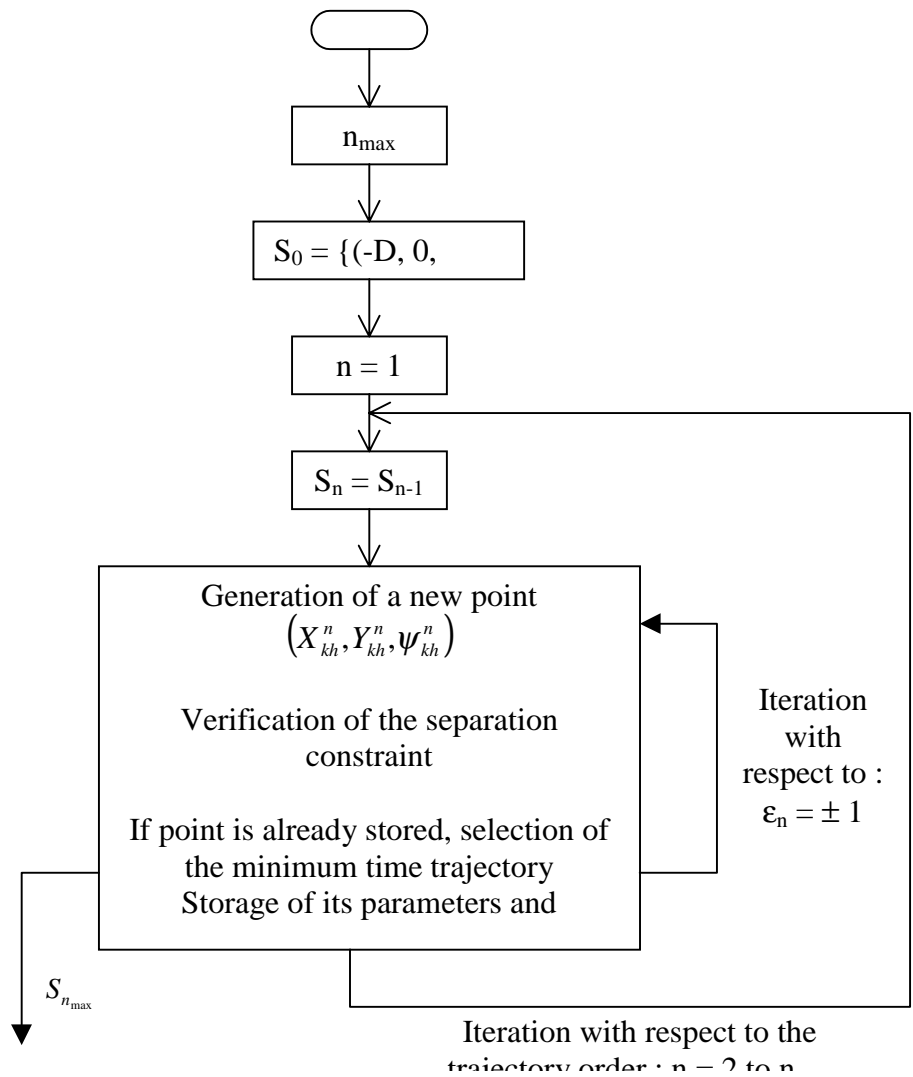

Figure 3. Reverse generation of regular convergence trajectories.

Then appears the problem of generation of training data for a neural structure. This point is tackled here using a reverse Dynamic Programming approach.

If a problem identical to the minimum time regular trajectory optimization problem, except that initial relative conditions are not prefixed (relation 18), is considered, it is possible to use in a simple way Dynamic Programming in the reverse direction to generate in a sequential manner a set of points $(x_0^P, y_0^P, \psi_0^P)$ which start at the final point of convergence and are reached through a n order regular trajectory.

Here possible headings and lengths of the straight line segments are discretized ($\theta_k = (k-1)\Delta\theta$, k = 1 to N and $l_h = (h-1)\Delta l$, h = 1 to M). The structure of the proposed iterative algorithm to generate such points is given in figure 3.

The set of generated trajectories is represented by a tree in $\Re^3$ with root at (-D, 0, $\psi_L$).

The next figure provides a view of the convergence area obtained by reverse Dynamic Programming.

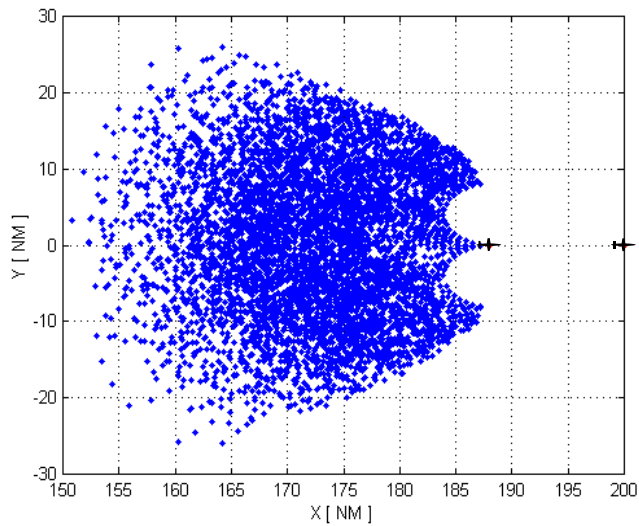

**Figure 4. Representation of the convergence area.**

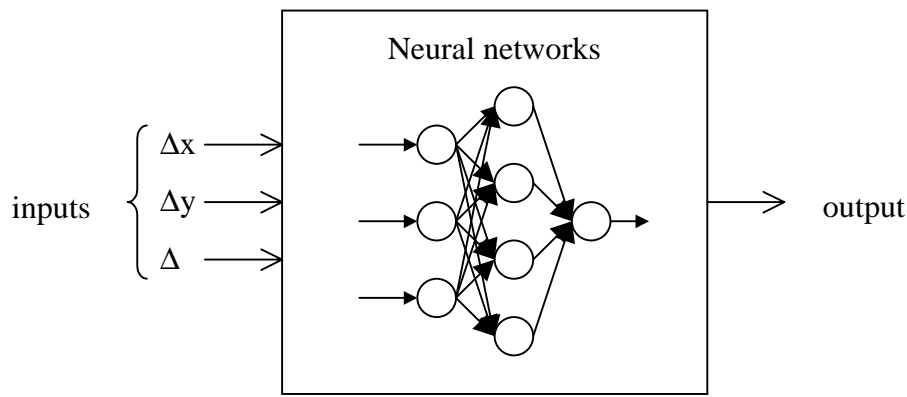

**Figure 5. The feedforward structure of the adopted neural networks.**

So, a data base can be made available for the training of a neural structure which will provide a solution to the original problem each time a convergence situation is proposed.

The reverse Dynamic Programming generation of convergence trajectories provides input/output pairs composed of :

Outputs : the relative initial position and heading of the Pursuer.

Inputs : the parameters of the minimum time convergence trajectory.

The general structure of the neural networks adopted here is of the feedforward class [16] and is shown in figure 5.

## 5   Simulation Results

A simulation study has been performed considering two Airbus A300 aircraft. It has been supposed that after the beginning of the convergence maneuver, the Leader aircraft makes a right turn to take a new constant heading. The guidance system of the Pursuer makes use of a neural trajectory generator to define, every second, new references (either a turn rate or a constant heading) for the autopilot. The guidance laws implemented in its autopilot are classical superposed PID loops (a fast piloting loop and a slower guidance loop) similar to those encountered in aircraft of this class [17].

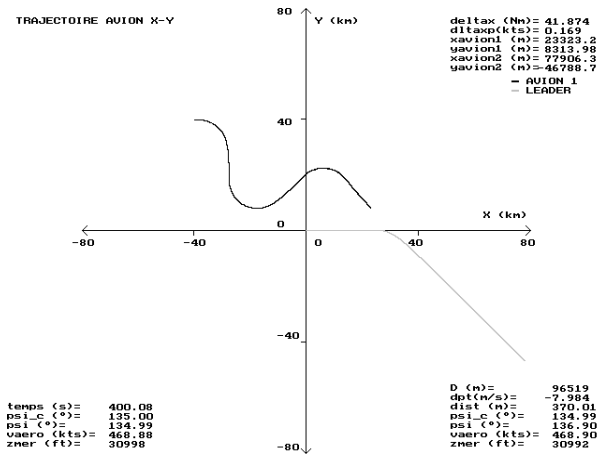

Fig. 6a. Simulation without Leader's intent information.

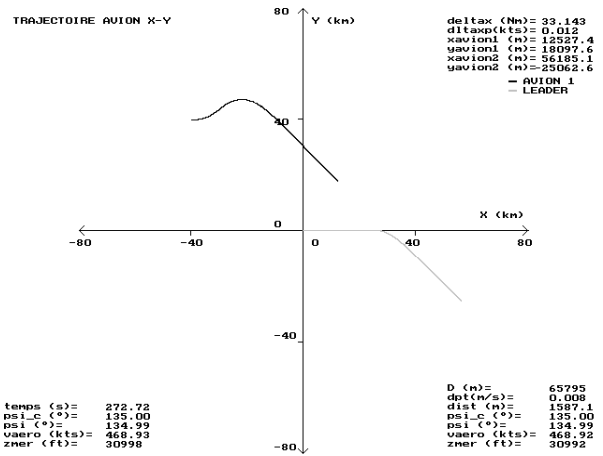

Fig. 6b. Simulation with Leader's intent information.

Simulation results are shown in Fig. 6a and 6b. Initial position of the Leader is taken as coordinate (0, 0) with a heading of $90^o$. The Pursuer is initially at 40 km behind and 40 km North of the Leader, with an absolute heading of $90^o$. In the case of Fig. 6a, it is assumed that the Pursuer has no knowledge of the Leader's flight plan. With the activation of the relative guidance mode, the Pursuer makes an initial right turn and starts convergence. When the Leader modifies its heading to $135^o$, the Pursuer changes progressively its own heading and then after convergence turns right to follow the Leader. This results in a very swaying trajectory, which is much undesirable in ATC standards. While in the case of Fig. 6b, the Pursuer knows exactly the intent of the Leader so, the pursuit trajectory can be generated taking into account the final route of the Leader. In the second case, the Pursuer avoids excessive maneuvering while the convergence time is much smaller than in the previous case.

## 6 Perspective

The proposed approach allows the development of a new on board system whose function can be :
- to provide directives to an on board guidance system such as a Flight Director operating in a "Relative Guidance" mode when the aircraft is under manual control,
- to provide references for the autopilot operating in a "Relative Guidance" mode when the aircraft is under its control. In this case, this function is very similar with other functions already available on advanced flight management systems.
- to visualize the convergence trajectory in real time on a navigation screen, with a display such as the one suggested in figure 7.

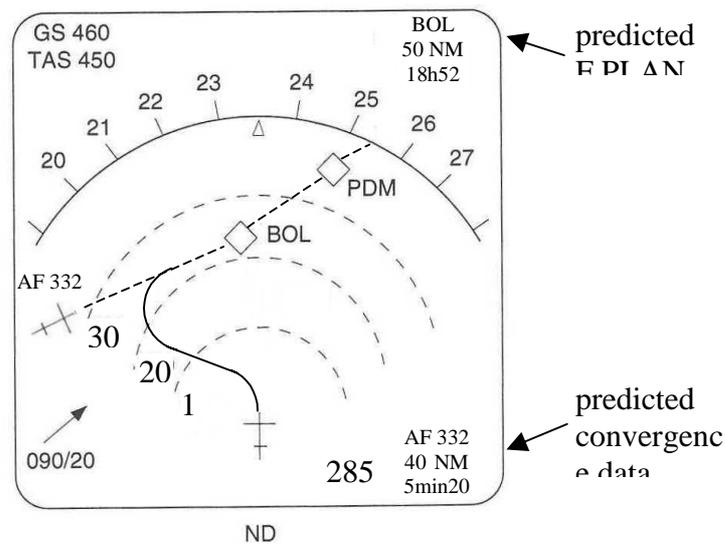

Figure 7. Presentation of a proposed or planned convergence maneuver to the pilot.

- to transmit the intentions of the Pursuer (convergence trajectory) to the Air Traffic Control system so that it can validate the convergence maneuver and integrate it into the traffic.